\documentclass[12pt]{article}
\usepackage[]{amsmath,amssymb}
\usepackage{amscd}
\usepackage{latexsym}
\usepackage{cite}

\newtheorem{definition}{Definition}[section]
\newtheorem{theorem}[definition]{Theorem}
\newtheorem{lemma}[definition]{Lemma}
\newtheorem{corollary}[definition]{Corollary}

\newtheorem{conjecture}[definition]{Conjecture}
\newtheorem{problem}[definition]{Problem}
\newtheorem{note}[definition]{Note}

\newtheorem{notation}[definition]{Notation}

\def\I{\mathbb I}
\def\N{\mathbb N}

\def\C{\mathbb C}

\def\Z{\mathbb Z}
\def\F{\mathbb F}

\def\T{\mathbb T}

\def\A{\mathbb A}
\def\B{\mathbb B}

\begin{document}
\title{\bf  
Augmented down-up algebras \\
and 
uniform posets
}
\author{
Paul Terwilliger and
Chalermpong Worawannotai
}
\date{}

\maketitle
\begin{abstract}
Motivated by the structure of the uniform posets we
introduce the notion of an augmented down-up (or ADU) algebra.
We discuss how ADU algebras are related to the
down-up algebras defined by Benkart and Roby.
For each ADU algebra
we give two presentations by generators and relations.
We also display a $\Z$-grading  and a linear basis.
In addition
we show that the center is
 isomorphic to a polynomial algebra in two variables.
We display  seven families of uniform posets and
show that each 
gives an ADU algebra module in a natural way.
The main inspiration for the ADU algebra concept
comes from the second author's thesis 
concerning a type
of uniform poset constructed using a dual polar graph.



 \end{abstract}

\section{Introduction}

\noindent In \cite{uniform} the first author introduced
the notion of a uniform poset, and constructed eleven families
of examples
from the classical
geometries. Among the examples are the polar spaces 
${\rm Polar}_b(N,\epsilon)$ and  the attenuated spaces
${\rm A}_b(N,M)$, as well as the posets 
${\rm Alt}_b(N)$,
${\rm Her}_q(N)$, and
${\rm Quad}_b(N)$ associated with the alternating, Hermitean, and
quadratic forms.
Another example is Hemmeter's poset
${\rm Hem}_b(N)$.
In 
\cite[Lemma~26.4]{boyd} the second author constructed a new family
of uniform posets using the dual polar graphs.
We denote these posets by
 ${\rm Polar}^{\rm top}_b(N,\epsilon)$ and describe
 them in
   Section 5 below.

\medskip
\noindent 
In \cite{benkart} Benkart and Roby introduced
the down-up algebras, and obtained
modules for these algebras
using
${\rm Alt}_b(N)$,
${\rm Her}_q(N)$,
${\rm Quad}_b(N)$, and
${\rm Hem}_b(N)$.
A 
down-up algebra module 
 is  obtained from ${\rm Polar}^{\rm top}_b(N,\epsilon)$
in a similar way. 
However, 
it appears that the down-up algebra
concept 
is not sufficiently robust to handle 
${\rm Polar}_b(N,\epsilon)$ or
${\rm A}_b(N,M)$. The same can be said for the generalized
down-up algebras
\cite{cas} 
and all the related algebras
\cite{joseph,
spsmith,
benkartkang,
cast, 
tang,
bavula,
bavoy,
delb,
fairlie,
flat,
fuj,
havlicek,
bruyn,
bruyn2,
odesskii,
woro,
van,
witten1,
witten} that we are aware of.
In the present paper we introduce a family of
algebras called augmented down-up algebras,
or ADU algebras for short.
We show that each of the uniform posets 
${\rm Polar}_b(N,\epsilon)$, 
${\rm A}_b(N,M)$,
${\rm Alt}_b(N)$,
${\rm Her}_q(N)$,
${\rm Quad}_b(N)$,
${\rm Hem}_b(N)$, 
 ${\rm Polar}^{\rm top}_b(N,\epsilon)$ 
gives an ADU algebra module in a natural way.

\medskip
\noindent
The 
ADU algebras are related to the down-up algebras as follows.
Given scalars $\alpha,\beta,\gamma$ the corresponding down-up
algebra $A(\alpha,\beta,\gamma)$ is 
defined by generators $e,f$ and relations
\begin{eqnarray*}
&&e^2 f = \alpha efe + \beta f e^2 + \gamma e,
\\
&&ef^2 = \alpha fef + \beta f^2 e + \gamma f.
\end{eqnarray*}
See \cite[p.~308]{benkart}.
To turn this into an ADU algebra we make
three adjustments as follows.
Let  $q$ denote a nonzero scalar that is not a root of unity.
We first require
\begin{eqnarray*}
\alpha = q^{-2s}+q^{-2t},
\qquad \qquad
\beta= -q^{-2s-2t}
\end{eqnarray*}
where $s,t$ are distinct integers.
Secondly, we 
add two  generators 
$k^{\pm 1}$ such that
$ke=q^2ek$ and $kf=q^{-2}fk$. Finally
we reinterpret $\gamma$ as a Laurent polynomial in
$k$ for which the coefficients of $k^s, k^t$ are zero.

\medskip
\noindent
From the above description
the ADU algebras are reminiscent
of the  quantum univeral enveloping algebra
$U_q({\mathfrak{sl}_2})$. To illuminate
the difference between these algebras, 
consider their center.
By
\cite[p.~27]{jantzen}
the center of
$U_q({\mathfrak{sl}_2})$
 is isomorphic to
a polynomial algebra in one variable.
As we will see, the center of an ADU algebra is isomorphic
to a polynomial algebra in two variables.

\medskip
\noindent The results of the present paper
are summarized as follows.
We define two algebras by generators and relations,
and show that they are isomorphic. We call the common
resulting algebra an ADU algebra.
For each ADU algebra we display
a $\Z$-grading and a linear basis.
We also
show that the center is isomorphic to a polynomial
algebra in two variables. 
We obtain ADU algebra modules from each of the
above seven examples of uniform posets.

\medskip
\noindent 
Recall 
the
natural numbers $\N=\lbrace 0,1,2,\ldots \rbrace$
and
integers $\Z=\lbrace 0,\pm 1,\pm 2,\ldots\rbrace$.

\section{Augmented down-up  algebras}

Our conventions for the paper are as follows.
An algebra is meant to be associative and
have a 1.
A subalgebra has the same 1 as the parent algebra.
Let $\F$ denote a field.
Let $\lambda$ denote an indeterminate.
Let $\F\lbrack \lambda, \lambda^{-1}\rbrack$
denote the $\F$-algebra of Laurent polynomials in $\lambda$
that have all coefficients in $\F$.
Pick $
\psi\in 
\F\lbrack \lambda, \lambda^{-1}\rbrack$ and write
$\psi = \sum_{i\in \Z} \alpha_i \lambda^i$.
By the {\it support} of $\psi$ we mean
the set 
$\lbrace i \in \Z | \alpha_i \not=0\rbrace$.
This set 
is finite.


\medskip
\noindent 
 Fix distinct  $s,t \in \Z$.
Define
\begin{eqnarray*}
\F\lbrack \lambda, \lambda^{-1}\rbrack_{s,t}
= {\rm Span} \lbrace \lambda^i | i \in \Z, \;i\not=s, i\not=t\rbrace.
\end{eqnarray*}
Note that
\begin{eqnarray*}
\F\lbrack \lambda, \lambda^{-1}\rbrack = 
\F\lbrack \lambda, \lambda^{-1}\rbrack_{s,t} + \F \lambda^s + \F \lambda^t
\qquad 
\hbox{\rm (direct sum).}
\end{eqnarray*}


\noindent
For $\psi \in 
\F\lbrack \lambda, \lambda^{-1}\rbrack$
the following are equivalent:
(i) $\psi \in 
\F\lbrack \lambda, \lambda^{-1}\rbrack_{s,t}$;
(ii) 
the integers $s,t$ are not in the support of
$\psi$.

%
%

\medskip
\noindent 
Fix a nonzero $q \in \F$ that is not a root of unity.

\begin{definition} 
\label{def:v1}
\rm
For
$\varphi \in 
\F\lbrack \lambda, \lambda^{-1}\rbrack_{s,t}$ 
the $\F$-algebra $\A=\A_q(s,t,\varphi)$ has generators
$e,f,k^{\pm 1}$ and relations
\begin{eqnarray}
&&kk^{-1}=1, \qquad k^{-1}k=1,
\nonumber
\\
&&ke=q^2ek, \qquad kf=q^{-2}fk,
\nonumber
\\
&&
e^2 f
-
(q^{-2s}+q^{-2t})efe
+
q^{-2s-2t}fe^2
=
e\varphi(k),  
\label{eq:eef}
\\
&&
ef^2
-
(q^{-2s}+q^{-2t})fef
+
q^{-2s-2t}f^2e
= \varphi(k) f.
\label{eq:eff}
\end{eqnarray}
\end{definition}

\begin{note}\rm
Referring to Definition
\ref{def:v1}, consider the special case in which
$\varphi \in \F$. Then 
the relations 
(\ref{eq:eef}), (\ref{eq:eff})
become the defining relations for the down-up
algebra $A(q^{-2s}+q^{-2t},-q^{-2s-2t},\varphi)$.
\end{note}

\begin{definition}
\label{def:v2}
\rm
For $\phi \in 
\F\lbrack \lambda, \lambda^{-1}\rbrack_{s,t}$
the $\F$-algebra $\B=\B_q(s,t,\phi)$ has generators
$C_s$, $C_t$, $E$, $F$, $K^{\pm 1}$ and relations
\begin{eqnarray}
&& C_s, \, C_t \quad {\hbox{\rm are central}},
\nonumber
\\
&&KK^{-1}=1, \qquad K^{-1}K=1,
\nonumber
\\
&&KE=q^2EK, \qquad KF=q^{-2}FK,
\nonumber
\\
&&FE = C_s q^s K^s + 
C_t q^t K^t
+ 
\phi(q K),
\label{eq:FE}
\\
&&EF = C_s q^{-s}K^s + C_t q^{-t} K^t + 
\phi(q^{-1}K).
\label{eq:EF}
\end{eqnarray}
\end{definition}

\noindent 
Next we describe how the algebras
in Definition
\ref{def:v1} and
 Definition
\ref{def:v2} are related. 

\begin{definition}
\label{def:phist}
\rm
We define an $\F$-linear map
$\F\lbrack \lambda, \lambda^{-1}\rbrack \to
\F\lbrack \lambda, \lambda^{-1}\rbrack$, $\psi\mapsto \psi_{s,t}$
as follows.
For $\psi \in 
\F\lbrack \lambda, \lambda^{-1}\rbrack$,
\begin{eqnarray*}
   \psi_{s,t}(\lambda) = \psi(q^{-1}\lambda)
   -(q^{-2s}+q^{-2t})\psi(q\lambda)+ q^{-2s-2t}\psi(q^3\lambda).
\end{eqnarray*}
\end{definition}

\noindent 
Recall the basis
$\lbrace \lambda^i \rbrace_{i \in \Z}$ for
$\F\lbrack \lambda, \lambda^{-1}\rbrack$.

\begin{lemma} 
\label{lem:stmap}
Consider the map
$\psi \mapsto \psi_{s,t}$ 
from Definition
\ref{def:phist}. For $i \in \Z$
the vector $\lambda^i$ is an eigenvector for the map.
The corresponding
eigenvalue is
$q^{3i}
(q^{-2i} - q^{-2s})
(q^{-2i} - q^{-2t})$. This eigenvalue is zero if and only if
$i\in \lbrace s,t\rbrace$.
\end{lemma}
\noindent {\it Proof:}
Use Definition
\ref{def:phist}.
\hfill $\Box$ \\

\noindent The following two lemmas are routine consequences
of Lemma
\ref{lem:stmap}.

\begin{lemma} 
For the map $\psi\mapsto \psi_{s,t}$ from
Definition
\ref{def:phist}
the image  is
$\F\lbrack \lambda, \lambda^{-1}\rbrack_{s,t}$
and the kernel is 
$\F \lambda^s + \F \lambda^t$.
\end{lemma}

\begin{lemma} 
For the map $\psi\mapsto \psi_{s,t}$ from
Definition
\ref{def:phist}
the restriction to
$\F\lbrack \lambda, \lambda^{-1}\rbrack_{s,t}$
is invertible.
\end{lemma}

\noindent 
Let $\varphi, \phi \in 
\F\lbrack \lambda, \lambda^{-1}\rbrack_{s,t}$ such that
$\varphi= \phi_{s,t}$.
We are going to show
that the algebras
$\A_q(s,t,\varphi)$ and $\B_q(s,t,\phi)$ are isomorphic.

\begin{lemma}
\label{lem:CsCt} 
For $\phi \in 
\F\lbrack \lambda, \lambda^{-1}\rbrack_{s,t}$ the following
hold in 
$\B_q(s,t,\phi)$:
\begin{eqnarray}
&&C_s =
\frac{q^{-t}FE-q^t EF + 
q^t\phi(q^{-1}K)-q^{-t}\phi(qK)
}
{q^{s-t}-q^{t-s}} K^{-s},
\label{eq:Cs}
\\
&&
C_t =
\frac{q^{-s}FE-q^s EF + 
q^s \phi(q^{-1}K)-q^{-s}\phi(q K)
}
{q^{t-s}-q^{s-t}} K^{-t}.
\label{eq:Ct}
\end{eqnarray}
Moreover the algebra
$\B_q(s,t,\phi)$ is generated by $E,F,K^{\pm 1}$.
\end{lemma}
\noindent {\it Proof:}
We first verify 
(\ref{eq:Cs}). 
In the expression on the right in 
(\ref{eq:Cs}),
eliminate  $FE$ and $EF$ 
using 
(\ref{eq:FE}) and 
(\ref{eq:EF}).
 After a routine simplification 
(\ref{eq:Cs}) is verified.
The equation 
(\ref{eq:Ct}) is similarly verified. The last assertion
follows from 
(\ref{eq:Cs}),
(\ref{eq:Ct}).
\hfill $\Box$ \\

\begin{lemma}
\label{lem:CsCtII} 
For $\phi \in 
\F\lbrack \lambda, \lambda^{-1}\rbrack_{s,t}$ the following
hold in 
$\B_q(s,t,\phi)$:
\begin{eqnarray}
&&
E^2 F 
-
(q^{-2s}+q^{-2t})EFE
+
q^{-2s-2t}FE^2
=
E\varphi(K),
\label{eq:unif1}
\\
&&
EF^2
-
(q^{-2s}+q^{-2t})FEF
+
q^{-2s-2t}F^2E
= \varphi(K) F.
\label{eq:unif2}
\end{eqnarray}
In the above lines $\varphi = \phi_{s,t}$.
\end{lemma}
\noindent {\it Proof:}
We first verify 
(\ref{eq:unif1}). In the expression on the left in 
(\ref{eq:unif1}),
 view $E^2F=E(EF)$, $EFE=E(FE)$,
$FE^2=(FE)E$ and eliminate each parenthetical
expression 
using
(\ref{eq:FE}) and 
(\ref{eq:EF}).
Simplify the result using $KE=q^2EK$ along with
$\varphi=\phi_{s,t}$ 
and
Definition
\ref{def:phist}.
The equation
(\ref{eq:unif1}) is now verified.
The equation 
(\ref{eq:unif2}) is similarly verified.
\hfill $\Box$ \\

\noindent The following definition is motivated by
Lemma
\ref{lem:CsCt}.

\begin{definition} 
\label{def:csct}
\rm 
For $\varphi \in 
\F\lbrack \lambda, \lambda^{-1}\rbrack_{s,t}$ let
$c_s, c_t$ denote the following 
elements in 
$\A_q(s,t,\varphi)$:
\begin{eqnarray}
&&c_s  =
\frac{q^{-t}fe-q^t ef + 
q^t\phi(q^{-1}k)-q^{-t}\phi(qk)
}
{q^{s-t}-q^{t-s}} k^{-s},
\label{def:cs}
\\
&&
c_t  =
\frac{q^{-s}fe-q^s ef + 
q^s \phi(q^{-1}k)-q^{-s}\phi(q k)
}
{q^{t-s}-q^{s-t}} k^{-t}.
\label{def:ct}
\end{eqnarray}
In the above lines $\phi$
denotes the unique element in 
$\F\lbrack \lambda, \lambda^{-1}\rbrack_{s,t}$ such that
$\varphi=\phi_{s,t}$.
\end{definition}

\begin{lemma} 
\label{lem:csct} 
With the notation and assumptions of
Definition
\ref{def:csct}, the elements $c_s, c_t$ are central
in 
$\A_q(s,t,\varphi)$. Moreover 
\begin{eqnarray}
&&fe = c_s q^s k^s + 
c_t q^t k^t
+ 
\phi(q k),
\label{eq:fe}
\\
&&ef = c_s q^{-s}k^s + c_t q^{-t} k^t + 
\phi(q^{-1}k).
\label{eq:ef}
\end{eqnarray}
\end{lemma}
\noindent {\it Proof:}
We first show that $c_s$ is central 
in $\A_q(s,t,\varphi)$. 
To do this 
we show  $c_s e= ec_s$,
$c_s f= f c_s$,
 $c_s k= kc_s$. To verify these equations,
eliminate each 
occurrence of $c_s$
using
(\ref{def:cs}), and simplify the result using
the relations in
Definition \ref{def:v1}. We have shown
that $c_s$ is central in
 $\A_q(s,t,\varphi)$.  One similarly shows that
$c_t$ is central in
 $\A_q(s,t,\varphi)$.  
We now verify (\ref{eq:fe}). In the expression on
the right in 
 (\ref{eq:fe}),
eliminate $c_s, c_t$ using
(\ref{def:cs}),
(\ref{def:ct}). After a routine simplification
(\ref{eq:fe}) is verified.
The equation
(\ref{eq:ef}) is similarly verified.
\hfill $\Box$ \\

\begin{theorem} 
Given  
$\varphi, \phi \in 
\F\lbrack \lambda, \lambda^{-1}\rbrack_{s,t}$ such that
 $\varphi= \phi_{s,t}$.
Then there exists an 
$\F$-algebra isomorphism
$\A_q(s,t,\varphi) \to \B_q(s,t,\phi)$ that sends
\begin{eqnarray*}
e \mapsto E, 
\quad \qquad 
f \mapsto F, 
\quad \qquad 
k^{\pm 1} \mapsto K^{\pm 1}.
\end{eqnarray*}
The inverse isomorphism sends
\begin{eqnarray*}
&&C_s \mapsto  c_s,
\qquad \quad 
C_t \mapsto 
c_t,
\qquad \quad 
E \mapsto e, 
\qquad \quad 
F \mapsto f, 
\qquad \quad 
K^{\pm 1} \mapsto k^{\pm 1}
\end{eqnarray*}
where $c_s,c_t$ are from
Definition
\ref{def:csct}.
\end{theorem}
\noindent {\it Proof:}
Combine Lemmas
\ref{lem:CsCt}, 
\ref{lem:CsCtII}, 
\ref{lem:csct}. 
\hfill $\Box$ \\

\begin{definition}\rm
By an {\it augmented down-up algebra}
we mean an
algebra
$\A_q(s,t,\varphi)$ from 
Definition
\ref{def:v1} or an algebra
$\B_q(s,t,\phi)$ from 
Definition
\ref{def:v2}.
\end{definition}

\noindent
Consider the algebra $\B=\B_q(s,t,\phi)$ from Definition
\ref{def:v2}. In Section 3 
we are going to show that the elements
$C_s, C_t$ generate the center $Z(\B)$, and
that 
$Z(\B)$ is isomorphic to a polynomial algebra in
two variables.
Because of this and following
\cite[p.~27]{jantzen}, it seems appropriate to call $C_s, C_t$
the {\it Casimir elements} for $\B_q(s,t,\phi)$.

\section{A $\Z$-grading and  linear basis for $\B_q(s,t,\phi)$}

\noindent Recall the algebra 
$\B=\B_q(s,t,\phi)$ from Definition
\ref{def:v2}. In this section we display a
$\Z$-grading for $\B$. We also display
a basis for the $\F$-vector space
 $\B$.

\medskip
\noindent 
Let $\mathcal A$ denote an
$\F$-algebra.
By a {\it $\Z$-grading of $\mathcal A$} we mean a sequence
 $\lbrace {\mathcal A}_n \rbrace_{n\in \Z}$ consisting of
  subspaces of $\mathcal A$
   such
   that
   \begin{eqnarray*}
   {\mathcal A} = \sum_{n\in \Z} \mathcal A_n
   \qquad \qquad
   {\mbox {\rm {(direct sum)}}},
   \end{eqnarray*}
   and $\mathcal A_m \mathcal A_n \subseteq \mathcal A_{m+n}$
   for all $m,n \in \Z$. Let
   $\lbrace {\mathcal A}_n \rbrace_{n\in \Z}$ denote
   a $\Z$-grading of $\mathcal A$.
   For $n \in \Z$
   we call
   $\mathcal A_n$ the {\it $n$-homogeneous component} of
   $\mathcal A$.
   We refer to $n$ as the {\it degree}
   of $\mathcal A_n$. An element of $\mathcal A$ is said to be
   {\it homogeneous of degree $n$} whenever it is contained in
   $\mathcal A_n$.

\begin{theorem}
\label{thm:basisg}
The 
algebra $\B$ 
has a $\Z$-grading
$\lbrace \B_n\rbrace_{n\in \Z}$
with the following properties:
\begin{enumerate}
\item[\rm (i)] 
The $\F$-vector space $\B_0$ has a basis
\begin{eqnarray}
K^h C_s^i C_t^j  \qquad    h\in \Z, \quad i,j \in \N.
\label{eq:zero}
\end{eqnarray}
\item[\rm (ii)] 
For $n\geq 1$, the $\F$-vector space $\B_n$ has a basis
\begin{eqnarray}
F^n K^h C_s^i C_t^j  \qquad    h\in \Z, \quad i,j \in \N.
\label{eq:pos}
\end{eqnarray}
\item[\rm (iii)] 
For $n\geq 1$, the $\F$-vector space $\B_{-n}$ has a basis
\begin{eqnarray}
E^n K^h C_s^i C_t^j  \qquad    h\in \Z, \quad i,j \in \N.
\label{eq:neg}
\end{eqnarray}
\end{enumerate}
Moreover the union of
{\rm (\ref{eq:zero})--(\ref{eq:neg})} is a basis for
the $\F$-vector space $\B$.
\end{theorem}
\noindent {\it Proof:}
Routinely applying the Bergman diamond lemma
\cite[Theorem~1.2]{berg}
one finds that the union of
(\ref{eq:zero})--(\ref{eq:neg}) is a basis for
the $\F$-vector space $\B$.
Let $\B_0$ denote the subspace of 
$\B$ spanned by 
(\ref{eq:zero}).
For $n\geq 1$ let $\B_n$ and $\B_{-n}$ denote
the subspaces of $\B$ spanned by
(\ref{eq:pos}) and
(\ref{eq:neg}), respectively. 
We show that
$\lbrace \B_n\rbrace_{n\in \Z}$ is a 
$\Z$-grading of $\B$.
By construction the sum
$\B = \sum_{n \in \Z} \B_n$ is direct.
By construction and since
$C_s, C_t$ are central we have
$C_s \B_n \subseteq \B_n$
and
$C_t \B_n \subseteq \B_n$
for $n \in \Z$.
Using $KE=q^2EK$ and
 $KF=q^{-2}FK$ we find
$K^{\pm 1} \B_n \subseteq \B_n$ for $n \in \Z$.
Using 
(\ref{eq:FE}) and 
(\ref{eq:EF})
we find
$E \B_n \subseteq \B_{n-1}$ and
$F \B_n \subseteq \B_{n+1}$ for
$n\in \Z$. By these comments and the construction
we see that $\B_m \B_n \subseteq \B_{m+n}$
for all $m,n\in\Z$. Therefore
$\lbrace \B_n\rbrace_{\in \Z}$ is a $\Z$-grading
of $\B$. The result follows.
\hfill $\Box$ \\

\noindent We emphasize a few points from Theorem
\ref{thm:basisg}.

\begin{corollary} With respect to the above
$\Z$-grading of $\B$,
the generators $C_s,C_t,E,F,K^{\pm 1}$ are homogeneous with
the following degrees:

\begin{center}
\begin{tabular}{c|c c c c c}
$v$ & $C_s$ &  $C_t$ & $E$  & $F$
& 
$K^{\pm 1}$ 
\\
\hline
{\rm degree of $v$} & $0$ & $0$ & $-1$ & $1$ & $0$
\end{tabular}
 \end{center}

\end{corollary}

\begin{corollary} The homogeneous component $\B_0$ is
the subalgebra of $\B$ generated by $C_s, C_t, K^{\pm 1}$.
The algebra $\B_0$ is commutative.
\end{corollary}

\noindent Let $\lbrace \lambda_i \rbrace_{i=0}^2$ denote
mutually commuting indeterminates.

\begin{corollary} There exists an $\F$-algebra isomorphism
$\B_0 \to \F\lbrack \lambda_0^{\pm 1},\lambda_1, \lambda_2\rbrack$
that sends
\begin{eqnarray*}
K^{\pm 1} \mapsto \lambda^{\pm 1}_0,
\qquad \quad 
C_s \mapsto \lambda_1,
\qquad \quad 
C_t \mapsto \lambda_2.
\end{eqnarray*}
\end{corollary}

\noindent The 
$\Z$-grading $\lbrace \B_n\rbrace_{n\in \Z}$ has the
following interpretation.

\begin{lemma}
\label{lem:bnmeaning}
Consider the $\F$-linear map
$\B \to \B$, $\xi \mapsto K^{-1}\xi K$.
For $n \in \Z$ the $n$-homogeneous component
$\B_n$ is an eigenspace of this map. The
corresponding eigenvalue is $q^{2n}$.
\end{lemma}
\noindent {\it Proof:}
Use the basis for $\B_n$ given in Theorem
\ref{thm:basisg}, along with the relations
$KE=q^2EK$ and $KF=q^{-2}FK$.
\hfill $\Box$ \\

\begin{corollary} 
\label{cor:B0meaning} 
The homogeneous component
$\B_0$ consists of the elements in $\B$ that
commute with $K$.
\end{corollary}
\noindent {\it Proof:}
Immediate from Lemma
\ref{lem:bnmeaning}.
\hfill $\Box$ \\

\section{The center of $\B_q(s,t,\phi)$}

\noindent Recall the algebra
$\B=\B_q(s,t,\phi)$ from Definition
\ref{def:v2}. In this section we describe
the center $Z(\B)$.

\begin{theorem} 
\label{thm:centerg}
The following is a basis for 
the $\F$-vector space $Z(\B)$:
\begin{eqnarray}
C_s^i C_t^j \qquad i,j \in \N.
\label{eq:center}
\end{eqnarray}
\end{theorem}
\noindent {\it Proof:}
By Theorem
\ref{thm:basisg} the elements
(\ref{eq:center}) are linearly independent
over $\F$, so they form a basis
for a subspace of $\B$ which we denote
by $Z'$. We show $Z'=Z(\B)$.
The elements $C_s, C_t$ are central
in $\B$ 
so $Z'\subseteq Z(\B)$.
To obtain the reverse inclusion, pick
$\xi \in Z(\B)$. The element $\xi$ commutes
with $K$, so $\xi \in \B_0$ by Corollary
\ref{cor:B0meaning}.
Recall the basis
(\ref{eq:zero}) 
for $\B_0$. Writing $\xi$ in this basis,
we find $\xi = \sum_{h\in \Z} K^h \xi_h $
where $\xi_h \in Z' $ for $h \in \Z$.
Using $KE=q^2EK$ and $\xi E=E\xi$  we obtain
$0 = E \sum_{h \in \Z} K^h \xi_h(q^{2h}-1)$.
Combining this with 
Theorem \ref{thm:basisg} we find
$\xi_h = 0 $ for all nonzero $h \in \Z$. Therefore
$\xi=\xi_0 \in Z'$.
We have shown
 $Z'= Z(\B)$ and the result follows.
\hfill $\Box$ \\

\begin{corollary}
There exists an $\F$-algebra isomorphism
$Z(\B) \to \F\lbrack \lambda_1, \lambda_2\rbrack$
that sends
\begin{eqnarray*}
C_s \mapsto \lambda_1, \qquad \qquad C_t \mapsto \lambda_2.
\end{eqnarray*}
\end{corollary}

\section{Uniform posets}

Recall the algebras $\A_q(s,t,\varphi)$ from 
Definition \ref{def:v1}. In this section we discuss 
how these algebras are related to the uniform posets
\cite{uniform}.

\medskip
\noindent 
Throughout this section we assume that $\F$ is the complex number
field $\C$.
Let $P$ denote a finite ranked poset with fibers
$\lbrace P_i\rbrace_{i=0}^N$ 
\cite[p.~194]{uniform}.
 Let $\C P$ denote the
vector space over $\C$ with basis $P$.
Let ${\rm End}(\C P)$ denote the $\C$-algebra consisting
of all $\C$-linear maps from $\C P$ to $\C P$.
We now define three elements in 
${\rm End}(\C P)$ called the lowering, raising, and $q$-rank
operators.
For $x \in P$,
the  lowering operator sends $x$ to the
sum of the elements in $P$ that are covered by $x$.
The  raising operator sends  $x$ to the
sum of the elements in $P$ that cover $x$.
The  $q$-rank operator sends $x$ to $q^{N-2i}x$
where $x \in P_i$.

\medskip
\noindent In 
\cite{uniform}
we introduced a class of finite ranked posets said to 
be {\it uniform}.  We refer the reader 
to that article for a detailed
 description of these posets. 
See also
\cite[p.~306]{benkart}
and
\cite{stefko,
qmatroid}.
In  
\cite[Section~3]{uniform}
we gave eleven examples of uniform posets. We are going to show
that six of these examples give an
 $\A_q(s,t,\varphi)$-module.
These six examples 
are listed in the first six rows of the table below.
The remaining row of the table contains an example
 ${\rm Polar}^{\rm top}_b(N,\epsilon)$ which is defined as follows.
Start with the poset 
${\rm Polar}_b(N,\epsilon)$ which we denote by $P$.
Using $P$ we define an undirected
graph $\Gamma$ as follows.
The vertex set of $\Gamma$ consists of the
top fiber $P_N$ of $P$.
Vertices $y,z \in P_N$ are adjacent in $\Gamma$ 
whenever they are distinct and cover a common element of $P$.
 The graph $\Gamma$ is often called a {\it dual
polar graph}
\cite[p.~274]{bcn},
\cite[Section~16]{boyd}.
 Fix a vertex $x\in P_N$.
Using $x$ we define a partial order $\leq$ on $P_N$
as follows. For $y,z \in P_N$ let $y \leq z$ whenever
$\partial(x,y)+\partial(y,z)= \partial(x,z)$, where
$\partial$ denotes path-length distance in $\Gamma$. 
We have turned $P_N$ into a poset. We
call this poset ${\rm Polar}^{\rm top}_b(N,\epsilon )$.
Using 
\cite[Lemma~26.4]{boyd} one checks that
${\rm Polar}^{\rm top}_b(N,\epsilon )$ is uniform.

\begin{theorem} In each row of the table below we give an
example of a uniform poset $P$.
For each example we display
integers $s<t$ and
a Laurent polynomial $\varphi \in 
\F\lbrack \lambda, \lambda^{-1}\rbrack_{s,t}$.
In each case the vector space
$\C P$ becomes an $\A_q(s,t,\varphi)$-module such that
the generator $e$ (resp. $f$) (resp. $k$) acts on 
$\C P$ as the lowering (resp. raising) (resp. $q$-rank) 
operator for $P$. For convenience, for each example we
display the element $\phi \in 
\F\lbrack \lambda, \lambda^{-1}\rbrack_{s,t}$ such that
$\varphi=\phi_{s,t}$.

\begin{center}
\begin{tabular}{c| c c c c}
{\rm example} &   $s$ & $t$  & $\varphi $ & $\phi $
\\
\hline
\\
${\rm Polar}_b(N,\epsilon )$ 
  & $0$ & $1$ & 
$-(q+q^{-1})(q^{2N+1+2\epsilon }\lambda^2+q^{N-3}\lambda^{-1})$ &
$-\frac{q^{2N+2\epsilon}\lambda^2+q^{N-1}\lambda^{-1}}{(q-q^{-1})^2}$
\\
\\
${\rm A}_b(N,M)$ & $-1$ & $0$ &
$-(q+q^{-1})q^{N+2M+1} \lambda$ & 
$- \frac{q^{N+2M-1}}{(q-q^{-1})^2} \lambda$
\\
\\
${\rm Alt}_b(N)$ & $-2$ & $-1$ &
$-(q+q^{-1})q^{2N+1}$ & 
$- \frac{q^{2N-2}}{(q-q^{-1})^2}$
\\
\\
${\rm Her}_q(N)$ & $-2$ & $-1$ &
$-(q+q^{-1})q^{2N+2}$ & 
$- \frac{q^{2N-1}}{(q-q^{-1})^2}$
\\
\\
${\rm Quad}_b(N)$ & $-2$ & $-1$ &
$-(q+q^{-1})q^{2N+3}$ & 
$- \frac{q^{2N}}{(q-q^{-1})^2}$
\\
\\
${\rm Hem}_b(N)$ & $-2$ & $-1$ &
$-(q+q^{-1})q^{2N+1}$ & 
$- \frac{q^{2N-2}}{(q-q^{-1})^2}$
\\
\\
${\rm Polar}^{\rm top}_b(N,\epsilon )$ & $-2$ & $-1$ &
$-(q+q^{-1})q^{2N+3+2\epsilon }$ & 
$- \frac{q^{2N+2\epsilon }}{(q-q^{-1})^2}$
\end{tabular}
 \end{center}
In the above table $b=q^2$.

\end{theorem}
\noindent {\it Proof:}
For each example except the last,
our assertions follow routinely from
\cite[Theorem~3.2]{uniform}.
For the last example
${\rm Polar}^{\rm top}_b(N,\epsilon )$ our assertions follow from
\cite[Theorem~1.10]{boyd}. Note that the parameter denoted
$\epsilon $ in 
\cite[Theorem~1.10]{boyd} is one more than the parameter denoted 
$\epsilon $ in
\cite[p.~201]{uniform}.
\hfill $\Box$ \\


\section{Acknowledgement} 
The main inspiration for the ADU algebra concept
comes from the second author's thesis \cite{boyd} 
concerning the
uniform poset
 ${\rm Polar}^{\rm top}_b(N,\epsilon)$.
To be more precise, it was his discovery of two
central elements that he called $C_1,C_2$ 
\cite[Section~28]{boyd} that suggested to us
how to define an ADU algebra.

\noindent Paul Terwilliger \hfil\break
\noindent Department of Mathematics \hfil\break
\noindent University of Wisconsin \hfil\break
\noindent 480 Lincoln Drive \hfil\break
\noindent Madison, WI 53706-1388 USA \hfil\break
\noindent email: {\tt terwilli@math.wisc.edu }\hfil\break

\medskip
\noindent Chalermpong Worawannotai \hfil\break
\noindent Department of Mathematics \hfil\break
\noindent University of Wisconsin \hfil\break
\noindent 480 Lincoln Drive \hfil\break
\noindent Madison, WI 53706-1388 USA \hfil\break
\noindent email: {\tt worawann@math.wisc.edu }\hfil\break

\newpage

\end{document}